# The number of solutions of $\phi(x) = m$

By Kevin Ford

## Abstract

An old conjecture of Sierpiński asserts that for every integer $k \geqslant 2$, there is a number $m$ for which the equation $\phi(x) = m$ has exactly $k$ solutions. Here $\phi$ is Euler's totient function. In 1961, Schinzel deduced this conjecture from his Hypothesis H. The purpose of this paper is to present an unconditional proof of Sierpiński's conjecture. The proof uses many results from sieve theory, in particular the famous theorem of Chen.

## 1. Introduction

Of fundamental importance in the theory of numbers is Euler's totient function $\phi(n)$. Two famous unsolved problems concern the possible values of the function $A(m)$, the number of solutions of $\phi(x) = m$, also called the multiplicity of $m$. Carmichael's Conjecture ([1],[2]) states that for every $m$, the equation $\phi(x) = m$ has either no solutions $x$ or at least two solutions. In other words, no totient can have multiplicity 1. Although this remains an open problem, very large lower bounds for a possible counterexample are relatively easy to obtain, the latest being $10^{10^{10}}$ ([10, Th. 6]). Recently, the author has also shown that either Carmichael's conjecture is true or a positive proportion of all values $m$ of Euler's function give $A(m) = 1$ ([10, Th. 2]). In the 1950's, Sierpiński conjectured that all multiplicities $k \geqslant 2$ are possible (see [9] and [16]), and in 1961 Schinzel [17] deduced this conjecture from his well-known Hypothesis H [18].

SCHINZEL'S HYPOTHESIS H. *Suppose $f_1(n), \ldots, f_k(n)$ are irreducible, integer-valued polynomials (for integral $n$) with positive leading coefficients. Also suppose that for every prime $q$, $q \nmid f_1(n) \cdots f_k(n)$ for some $n$. Then the numbers $f_1(n), \ldots, f_k(n)$ are simultaneously prime for infinitely many positive integers $n$.*

Hypothesis H is a generalization of Dickson's prime $k$-tuples conjecture [7], which covers the special case when each $f_i$ is linear. Schinzel's construction of a number with multiplicity $k$ requires that as many as $k$ distinct polynomials with degrees as large as $3k$ simultaneously represent primes. Using a different, iterative argument, the author ([10, Th. 9]) deduced Sierpiński's Conjecture from the Prime $k$-tuples Conjecture, the argument requiring only three linear polynomials to be simultaneously prime. The proof is quite simple, and we reproduce it here.

LEMMA 1.1. *Suppose $A(m) = k$ and $p$ is a prime satisfying*

(i) $p > 2m + 1$,

(ii) $2p + 1$ *and* $2mp + 1$ *are prime*,

(iii) $dp + 1$ *is composite for all* $d \mid 2m$ *except* $d = 2$ *and* $d = 2m$.

*Then $A(2mp) = k + 2$.*

*Proof.* Suppose $\phi(x_1) = \cdots = \phi(x_k) = m$ and $\phi(x) = 2mp$. Condition (i) implies $p \nmid x$, hence $p \mid (q-1)$ for some prime $q$ dividing $x$. Since $(q-1) \mid 2mp$, we have $q = dp + 1$ for some $d \mid 2m$. We have $q > 2p$, so $q^2 \nmid x$ and thus $\phi(x) = (q-1)\phi(x/q)$. By conditions (ii) and (iii), either $q = 2p + 1$ or $q = 2mp+1$. In the former case, $\phi(x/q) = m$, which has solutions $x = (2p+1)x_i$ ($1 \leq i \leq k$). In the latter case, $\phi(x/q) = 1$, which has solutions $x = q$ and $x = 2q$. □

Now suppose $A(m) = k$ and let $d_1, \ldots, d_j$ be the divisors of $2m$ with $3 \leq d_i < 2m$. Let $p_1, \ldots, p_j$ be distinct primes satisfying $p_i > d_i$ for each $i$. Using the Chinese Remainder Theorem, let $a \mod b$ denote the intersection of the residue classes $-d_i^{-1} \mod p_i$ ($1 \leq i \leq j$). Then for every $h$ and $i$, $(a + bh)d_i + 1$ is divisible by $p_i$, hence composite for large enough $h$. The Prime $k$-tuples Conjecture implies that there are infinitely many numbers $h$ so that $p = a + hb$, and $2p + 1$ and $2mp + 1$ are simultaneously prime. By Lemma 1.1, $A(2mp) = k + 2$ for each such $p$. Starting with $A(1) = 2$ and $A(2) = 3$, Sierpiński's Conjecture follows by induction on $k$.

Although Hypothesis H has not been proved in even the simplest case of two linear polynomials (generalized twin primes), sieve methods have shown the conclusion to hold if the numbers $f_1(n), \ldots, f_k(n)$ are allowed to be primes or "almost primes" (nonprimes with few prime factors). See [13] for specific statements. Recently, the author and S. Konyagin [11] employed the theory of almost primes to prove that $A(m) = k$ is solvable for every even $k$. Although falling short of a complete proof of Sierpiński's Conjecture, the method did provide a solution to the corresponding problem for the sum of divisors function $\sigma(n)$. It is proved ([11], Theorem 1) that for every $k \geq 0$, there is a number $m$ for which $\sigma(x) = m$ has exactly $k$ solutions, settling another conjecture of Sierpiński. Curiously, the barrier to solving $A(m) = k$ for odd $k$ by this method is the probable nonexistence of a number with $A(m) = 1$.

Combining a different application of almost primes with the iterative approach from [10], we now provide an unconditional proof of Sierpiński's Conjecture.

THEOREM 1. *For each integer $k \geqslant 2$, there is an integer $m$ with $A(m) = k$.*

Our proof of Theorem 1 is by induction on $k$. We start with $A(10) = 2$, $A(2) = 3$ and for the induction step, we show:

THEOREM 2. *Suppose $A(m) = k$ and $m$ is even. Then there is a number $m'$ such that $A(mm') = k + 2$.*

Combining Theorem 1 with Theorem 2 of [10] gives a strong result on the distribution of totients (values taken by $\phi(n)$) with a given multiplicity.

COROLLARY 3. *For each $k \geqslant 2$, a positive proportion of all totients have multiplicity $k$. Specifically, if $V(x)$ is the number of totients $\leqslant x$, and $V_k(x)$ is the number of totients $\leqslant x$ with multiplicity $k$, then for each $k \geqslant 2$,*

$$\liminf_{x \to \infty} \frac{V_k(x)}{V(x)} > 0.$$

At the center of the proof of Theorem 2 is a famous theorem of Chen ([3], [4]).

LEMMA 1.2 (Chen's theorem). *For each even natural number $m$ and $x \geqslant x_0(m)$, there are $\gg x/\log^2 x$ primes $s \in (x/2, x]$, $s \equiv 1 \pmod{m}$ such that $(s-1)/m$ has at most two prime factors, each of which exceeds $x^{1/10}$.*

*Remark.* A proof of Chen's Theorem also appears in Chapter 11 of [13]. What Chen actually proved is that there are $\gg x/\log^2 x$ primes $p \leqslant x$ with $p + 2$ having at most two prime factors. Trivial modifications to the proof give Lemma 1.2.

HYPOTHESIS $S(m)$. *For some constants $c > 0$ and $x_0$, if $x > x_0$ then the number of primes $s \in (x/2, x]$ for which $(s-1)/m$ is also prime is $\geqslant cx/\log^2 x$.*

We suppose $\phi(\alpha) = m$ has the $k$ solutions $\alpha_1, \ldots, \alpha_k$. The proof is then broken into two parts: (I) assuming Hypothesis $S(m)$ is true; (II) assuming Hypothesis $S(m)$ is false. For Part I we show that for many such primes $s$, $A(mq(q-1)) = k + 2$, where $q = (s-1)/m$ and $\phi(n) = mq(q-1)$ has the $k + 2$ solutions $n = q^2 \alpha_i$ ($i = 1, \ldots, k$), $n = q(mq + 1)$ and $n = 2q(mq + 1)$ (the "trivial" solutions).

For Part II, Chen's Theorem tells us that for an unbounded set of $x$, there are $\gg x/\log^2 x$ primes $s \in (x/2, x]$ with $s = 1 + mpq$, where $p$ and $q$ are primes and $x^{1/10} \leqslant p \leqslant q \leqslant x^{9/10}$. We show that for many such $s$, $A(mpq(p-1)(q-1)) = k+2$, the solutions of

$$\phi(n) = mpq(p-1)(q-1) \tag{1.1}$$

being only the "trivial" ones $n = p^2q^2\alpha_i$ ($i = 1, \ldots, k$), $n = pq(mpq+1)$ and $n = 2pq(mpq+1)$. Proving this is more difficult than Part I, but a simple heuristic indicates that such $s$ are likely. The basis for the heuristic (and the actual proof) is an extension of the classical result of Hardy and Ramanujan [14] that most integers $n \leqslant x$ have about $\log\log x$ prime factors. Erdős [8] combined their method with Brun's sieve to show that for most primes $p \leqslant x$, $p-1$ has about $\log\log x$ prime factors. One can further refine these ideas (e.g. Lemma 2.8 of [10] or Lemma 4.3 below) to show that for most primes $p \leqslant x$, the prime factors of $p-1$ are nicely distributed. This means that $p-1$ has about $\log\log b - \log\log a$ prime factors in the interval $[a,b]$ for $1 \leqslant a \leqslant b \leqslant x$. Call such primes "normal". If $p$ and $q$ are both normal primes and $x$ is large, then $w = p(p-1)q(q-1)m$ has about $2\log\log x$ prime factors. Fix $t$ and suppose $n$ has $t$ distinct prime factors $r_1, \ldots r_t$, different from $p$ and $q$. For simplicity suppose $n$ is square-free. Then

$$\phi(n) = (r_1 - 1)\cdots(r_t - 1) = w.$$

For a given factorization $w = d_1 \cdots d_t$, the probability that $d_1 + 1, \ldots, d_t + 1$ are all prime is about $((\log d_1) \cdots (\log d_t))^{-1}$. The critical case is $t = 2$. Since $p$ and $q$ are normal, $w$ has about $2\log\log y$ prime factors $\leqslant y$, and thus has about $(\log y)^{2\log 2}$ divisors $\leqslant y$. By partial summation,

$$\sum_{\substack{w = d_1 d_2 \\ d_i \geqslant 2}} \frac{1}{(\log d_1)(\log d_2)} \ll (\log x)^{2\log 2 - 2},$$

the right side representing the probability that (1.1) has a nontrivial solution for a given pair $(p, q)$. Similarly, it can be shown that

$$\sum_{\substack{w = d_1 \cdots d_t \\ d_i \geqslant 2}} \frac{1}{(\log d_1) \cdots (\log d_t)} \ll (\log x)^{2\log 2 - 2},$$

the bulk of the sum coming from taking $d_1, \ldots, d_{t-2}$ all small.

The biggest practical problem with these arguments is the fact that $\phi(n)$ is quadratic in $q$ for Part I and quadratic in both $p$ and $q$ for Part II. Thus, when certain sieve estimates are employed, the bounds have factors of the form

$$P(D, z) = \prod_{p \leqslant z} \left(1 - \frac{1}{p}\left(\frac{D}{p}\right)\right),$$

where $(\frac{D}{p})$ is the Legendre symbol. Heuristically, such products should be $O(1)$ on average, and showing this for a wide range of the parameters $D$ and $z$ requires strong results on the distribution of primes in arithmetic progressions (see (2.1) and Lemmas 2.6 – 2.8 and 4.2 below).

In the next section, various lemmas needed in the proof of both parts are presented. Section 3 is devoted to the proof in Part I. Section 4 contains further lemmas needed for Part II, and the argument counting nontrivial solutions of (1.1) with $(n, pq) > 1$. The last section deals with the hardest case, namely counting solutions of (1.1) with $(n, pq) = 1$.

The following notation will be used throughout. Let $P^+(n)$ denote the largest prime factor of $n$ and $P^-(n)$ denote the smallest prime factor of $n$. The functions $\omega(n)$ and $\Omega(n)$ count the number of distinct prime factors of $n$, respectively the number of prime factors counted according to multiplicity. Let $\Omega(n, U, T)$ denote the number of prime factors $p$ of $n$ such that $U \leqslant p < T$, counted according to multiplicity. Let $\left(\frac{a}{n}\right)$ denote the Legendre-Jacobi symbol and let $\tau(n)$ be the number of positive divisors of $n$. Constants implied by the Landau $O-$, Vinogradov $\ll -$ and $\gg -$ symbols may depend on $m$, but not on any other parameter unless indicated. The symbols $p, q$, and $\varpi$ will always denote primes. The symbol $m$ will always refer to the number in Theorem 2. The cardinality of a set $S$ will be denoted as $\#S$ or $|S|$.

## 2. Preliminary lemmas

Aside from Chen's Theorem, the other sieve results we require all follow from a basic sieve inequality called the "Fundamental Lemma" of the combinatorial sieve (Theorem 2.5 of [13]).

LEMMA 2.1. *Suppose $F$ is a polynomial of degree $h \geqslant 1$ with integer coefficients and positive leading coefficient. Denote by $\varrho(\varpi)$ the number of solutions of $F(n) \equiv 0 \pmod{\varpi}$. For $2 \leqslant z \leqslant x$,*

$$\#\{n \leqslant x : P^-(F(n)) > z\} = x \prod_{\varpi \leqslant z} \left(1 - \frac{\varrho(\varpi)}{\varpi}\right)$$
$$\times \left\{1 + O_h\left(e^{-\sqrt{\log x}} + e^{-u(\log u - \log\log 3u - \log h - 2)}\right)\right\},$$

*where $u = \log x / \log z$.*

*Proof.* If $\varrho(\varpi) = \varpi$ for some prime $\varpi \leqslant z$, then the left side is zero. Otherwise the lemma follows from [13, Th. 2.6] and following remarks. □

*Remark.* The error term is $O_h(1)$, and this will suffice for most of our applications. Also, for future reference, we note that for $\varpi \nmid 2a$,

$$(2.1) \qquad \#\{n : an^2 + bn + c \equiv 0 \pmod{\varpi}\} = 1 + \left(\frac{b^2 - 4ac}{\varpi}\right).$$

LEMMA 2.2. *Suppose*

$$F(n) = (a_1 n + b_1) \cdots (a_h n + b_h),$$

*where the $a_i$ are positive integers and the $b_i$ are integers. Suppose*

$$E := a_1 \cdots a_h \prod_{1 \leqslant i < j \leqslant h} (a_i b_j - a_j b_i) \neq 0.$$

*If $2 \leqslant z \leqslant x$ and $\log E \leqslant \log^3 x$, then*

$$\#\{n \leqslant x : P^-(F(n)) \geqslant z\} \ll_h \frac{x}{(\log z)^h} \left(\frac{E}{\phi(E)}\right)^h$$
$$\ll_h x (\log \log x)^h (\log z)^{-h}.$$

*Proof.* In the notation of Lemma 2.1, we have $\varrho(\varpi) = h$ when $\varpi \nmid E$; therefore the number of $n$ in question is

$$\ll_h x \prod_{\substack{\varpi \leqslant z \\ \varpi \nmid E}} \left(1 - \frac{1}{\varpi}\right)^h \ll_h x \left(\frac{E}{\phi(E)}\right)^h (\log z)^{-h}.$$

The last part of the lemma follows from the inequality $E/\phi(E) \ll \log \log E$. $\square$

LEMMA 2.3. *Let $V(x)$ be the number of distinct values of $\phi(n)$ which are $\leqslant x$. Then for every $\varepsilon > 0$,*

$$V(x) \ll_\varepsilon \frac{x}{(\log x)^{1-\varepsilon}}.$$

*Proof.* This is the main theorem of [8]. $\square$

LEMMA 2.4. *The number of $n \leqslant x$ divisible by a number $s$ satisfying $s \geqslant \exp\{(\log \log x)^2\}$ and $P^+(s) \leqslant s^{10/\log \log x}$ is $\ll x/\log^{10} x$.*

*Proof.* For $\exp\{(\log \log x)^2\} \leqslant y \leqslant x$, let $G(y)$ be the number of $s \leqslant y$ with $P^+(s) \leqslant s^{10/\log \log x}$. Let $n$ denote a number with no prime factor exceeding

$W := y^{10/\log\log x}$, and set $\alpha = 1 - \frac{20\log\log x}{\log y}$. Then

$$G(y) \leqslant \sum_{n \leqslant y} \left(\frac{y}{n}\right)^\alpha < y^\alpha \sum_{n=1}^\infty n^{-\alpha} = \frac{y}{(\log x)^{20}} \prod_{p \leqslant W} (1-p^{-\alpha})^{-1}$$

$$\ll \frac{y}{(\log x)^{20}} \exp\left(\sum_{p \leqslant W} p^{-\alpha}\right).$$

By the prime number theorem,

$$\sum_{p \leqslant W} p^{-\alpha} \leqslant (1+o(1)) \int_2^W \frac{dt}{t^\alpha \log t} \leqslant (1+o(1)) \int_{(\log y)^{-1}}^{200} \frac{e^u\, du}{u}$$

$$\leqslant (1+o(1))(e^{200} + \log\log y) \leqslant 2\log\log x$$

for large $x$. Thus $G(y) \ll y/(\log x)^{18}$ and the lemma follows by partial summation. For a survey of results on the distribution of numbers without large prime factors, see [15]. □

The next four lemmas are needed to bound the product $P(D,z)$ mentioned in the introduction.

LEMMA 2.5. *When $(a,n)=1$, $n$ is odd and square-free,*

$$\left|\sum_{\substack{1 \leqslant j \leqslant n \\ (j,n)=1}} \left(\frac{aj+b}{n}\right)\right| \leqslant 1.$$

*When $(ab,n)=1$ and $n$ is square-free,*

$$\left|\sum_{1 \leqslant j \leqslant n} \left(\frac{j}{n}\right)\left(\frac{aj+b}{n}\right)\right| = 1.$$

*Proof.* Denote by $f(n;a,b)$ the first sum in the lemma. The desired bound follows from the identities $f(n;a,b) = -(\frac{b}{n})$ for prime $n$, and

$$f(n_1 n_2; a, b) = f(n_1; an_2, b)f(n_2; an_1, b)$$

when $(a, n_1 n_2) = (n_1, n_2) = 1$. The second assertion follows in a similar fashion. Here we adopt the convention that $(\frac{a}{1}) = 1$ for all $a$. □

LEMMA 2.6. *For some absolute constant $y_0$, if $y \geqslant y_0$, $1 \leqslant c \leqslant 2$, $E \geqslant 1$, and $D \in \mathbb{Z}$, then*

$$\prod_{\substack{\varpi \leqslant y \\ \varpi \nmid E}} \left(1 - \frac{1}{\varpi}\left(\frac{D}{\varpi}\right)\right)^c \leqslant 4 \sum_{\substack{n \leqslant y^{\log\log y} \\ (n,2E)=1 \\ P^+(n) \leqslant y}} \frac{c^{\omega(n)}\mu(n)}{n}\left(\frac{D}{n}\right).$$

*Proof.* If $D = 0$, the left side is 1 and the right side is 4, so suppose $D \neq 0$. For $1 \leqslant c \leqslant 2$ and $x \geqslant -1/3$, $(1+x)^c \leqslant (1+cx)(1+3x^2)$; thus

$$\prod_{\substack{\varpi \leqslant y \\ \varpi \nmid E}} \left(1 - \frac{1}{\varpi}\left(\frac{D}{\varpi}\right)\right)^c \leqslant \frac{9}{4} \prod_{\varpi \geqslant 3} \left(1 + \frac{3}{\varpi^2}\right) \prod_{\substack{3 \leqslant \varpi \leqslant y \\ \varpi \nmid E}} \left(1 - \frac{c}{\varpi}\left(\frac{D}{\varpi}\right)\right)$$

$$= \frac{9}{4} \prod_{\varpi \geqslant 3} \left(1 + \frac{3}{\varpi^2}\right) \sum_{\substack{P^+(n) \leqslant y \\ (n, 2E) = 1}} \frac{c^{\omega(n)} \mu(n)}{n} \left(\frac{D}{n}\right)$$

$$\leqslant 3.92 \sum_{\substack{P^+(n) \leqslant y \\ (n, 2E) = 1 \\ n \leqslant y^{\log \log y}}} \frac{c^{\omega(n)} \mu(n)}{n} \left(\frac{D}{n}\right) + 3.92 \sum_{\substack{n > y^{\log \log y} \\ P^+(n) \leqslant y}} \frac{c^{\omega(n)}}{n}.$$

Since $c^{\omega(n)} \leqslant \tau(n)$, by Lemma 2.4 and partial summation, we obtain

$$\sum_{\substack{n > y^{\log \log y} \\ P^+(n) \leqslant y}} \frac{c^{\omega(n)}}{n} \leqslant \sum_{\substack{d_1 d_2 > y^{\log \log y} \\ P^+(d_1 d_2) \leqslant y}} \frac{1}{d_1 d_2} \ll \sum_{\substack{d_2 > y^{(\log \log y)/2} \\ P^+(d_2) \leqslant y}} \frac{1}{d_2} \sum_{d_1 \leqslant d_2} \frac{1}{d_1}$$

$$\ll \sum_{\substack{d_2 > y^{(\log \log y)/2} \\ P^+(d_2) \leqslant d_2^{3/\log \log d_2}}} \frac{\log d_2}{d_2}$$

$$\ll (\log y)^{-8}.$$

As the left side in the lemma is $\gg (\log y)^{-c}$, the desired bound follows for sufficiently large $y$. In particular, the sum in the lemma is positive and $\gg (\log y)^{-c}$. □

LEMMA 2.7. *Suppose $y \geqslant y_0$, $U \geqslant y^{10(\log \log y)^2}$, $A, B$ and $R$ are positive integers with $R$ even and $\log(ABR) \leqslant \log^3 U$. Suppose $D(s) = A(A - Bs)$, $D(s) = As(As - B)$ or $D(s) = A \pm Bs$. Uniformly in $A, B, R$ and $1 \leqslant c \leqslant 2$,*

$$\sum_{\substack{s \leqslant U \\ s \text{ prime}}} \prod_{\varpi \leqslant y} \left(1 - \frac{1}{\varpi}\left(\frac{D(s)}{\varpi}\right)\right)^c \ll \frac{U(\log \log U)^{1+c}}{\log U},$$

$$\sum_{\substack{s \leqslant U \\ s, sR+1 \text{ prime}}} \prod_{\varpi \leqslant y} \left(1 - \frac{1}{\varpi}\left(\frac{D(s)}{\varpi}\right)\right)^c \ll \frac{U(\log \log U)^{4+c}}{(\log U)^2}.$$

*Proof.* Let $z = U^{1/\log\log U}$. If $D(s) = A(A - Bs)$ or $D(s) = As(As - B)$, set $E = AB$; otherwise set $E = B$. By imposing the addition condition that $\varpi \nmid E$, we see that the products on the left are only increased when $D(s) = A(A - Bs)$ or $D(s) = As(As - B)$, and in the case $D(s) = A + Bs$ the product is changed by a factor

$$\prod_{\varpi \mid B} \left(1 - \frac{1}{\varpi}\left(\frac{A}{\varpi}\right)\right)^c \ll (\log\log B)^c \ll (\log\log U)^c.$$

We will also replace the conditions that "$s$ prime" and "$s, sR + 1$ prime", with $P^-(s) > z$, respectively $P^-(s(sR+1)) > z$. This will eliminate at most $z$ values of $s$ from the sums on the left, reducing the sums by at most $O(z(\log y)^2) = O_\varepsilon(U^\varepsilon)$. Let $\mathscr{S}_1$ denote the set of $s \leq U$ with $P^-(s) > z$ and let $\mathscr{S}_2$ denote the $s \leq U$ with $P^-(s(sR+1)) > z$. By Lemma 2.6,

$$\sum_{s\in\mathscr{S}_1} \prod_{\substack{\varpi \leq y \\ \varpi \nmid E}} \left(1 - \frac{1}{\varpi}\left(\frac{D(s)}{\varpi}\right)\right)^c \leq 4 \sum_{s\in\mathscr{S}_1} \sum_{\substack{n \leq y^{\log\log y} \\ (n,2E)=1 \\ P^+(n) \leq y}} \frac{c^{\omega(n)}\mu(n)}{n}\left(\frac{D(s)}{n}\right)$$

$$\leq 4 \sum_{\substack{n \leq y^{\log\log y} \\ (n,2E)=1}} \frac{c^{\omega(n)}\mu^2(n)}{n}\left|\sum_{s\in\mathscr{S}_1}\left(\frac{D(s)}{n}\right)\right|.$$

Denote by $S(\mathscr{S}; n, b)$ the number of elements of the set $\mathscr{S}$ which lie in the residue class $b \bmod n$. By the sieve fundamental lemma (Lemma 2.1), and the fact that $n \leq z^{1/5}$, uniformly in $(b, n) = 1$ we have

$$S(\mathscr{S}_1; n, b) = \frac{UW_1}{n}\left(1 + O((\log U)^{-10})\right),$$

where

$$W_1 = \prod_{\substack{\varpi \leq z \\ \varpi \nmid n}} (1 - 1/\varpi) \ll \frac{n}{\phi(n)}\frac{\log\log U}{\log U}.$$

For each $n$, Lemma 2.5 gives

$$\sum_{s\in\mathscr{S}_1}\left(\frac{D(s)}{n}\right) = \sum_{(b,n)=1} S(\mathscr{S}_1; n, b)\left(\frac{D(b)}{n}\right)$$

$$= \frac{UW_1}{n}\sum_{(b,n)=1}\left(\frac{D(b)}{n}\right) + O(U(\log U)^{-10})$$

$$\ll \frac{U\log\log U}{\log U}\frac{1}{\phi(n)} + \frac{U}{\log^{10} U},$$

and the first inequality follows. The second inequality follows in the same manner by summing over $s \in \mathscr{S}_2$ and using

$$S(\mathscr{S}_2; n, b) = \frac{UW_2}{n}(1 + O((\log U)^{-10})),$$

where $W_2 = 0$ if $(n, 1 + bR) > 1$ and otherwise

$$W_2 = \prod_{\substack{\varpi \leqslant z \\ \varpi \nmid nR}} \left(1 - \frac{2}{\varpi}\right) \prod_{\substack{\varpi \leqslant z \\ \varpi | R, \varpi \nmid n}} \left(1 - \frac{1}{\varpi}\right) \ll \left(\frac{nR}{\phi(nR)}\right)^2 \frac{(\log \log U)^2}{\log^2 U}. \quad \square$$

LEMMA 2.8. *For some absolute constant $b > 0$, the following holds. For each $Q \geqslant 3$, there is a number $d_Q \gg_A (\log Q)^A$ for every $A > 0$, so that whenever $|D| \leqslant Q$, $d_Q \nmid 8D$, and $\log z \leqslant \exp\{(\log Q)^b\}$,*

$$\prod_{Q^{\log \log Q} < \varpi \leqslant z} \left(1 - \frac{1}{\varpi}\left(\frac{D}{\varpi}\right)\right) \ll 1.$$

*Furthermore, $16 \nmid d_Q$ and $p^2 \nmid d_Q$ for all odd primes $p$.*

*Remark.* From now on $d_Q$ will refer to the number in this lemma. The implied constant in $d_Q \gg_A (\log Q)^A$ is ineffective for $A \geqslant 2$. With more work, the set of exceptional $D$ may be reduced to those with $f(D) = f(d_Q)$, where $f(n)$ is the product of the odd primes which divide $n$ to an odd power. The upper limit of $z$ may also be extended easily.

*Proof.* We may assume without loss of generality that $D$ is square-free. If $|D| \leqslant 2$ the result follows from the prime number theorem for arithmetic progressions, since $\left(\frac{D}{\varpi}\right)$ then depends only on $\varpi$ modulo 8. Next suppose $|D| \geqslant 3$. The desired bound follows from a version of the prime number theorem for arithmetic progressions deducible from the Landau-Page theorem and a zero density result of Gallagher (see [6] and [12]). The number $d_Q$ is the order of a real, primitive Dirichlet character whose corresponding $L$-function has a real zero $\beta > 1 - c_0/\log Q$. If $c_0 > 0$ is a sufficiently small absolute constant, the Landau-Page theorem asserts that there is at most one such primitive character of order $\leqslant 8Q$. Furthermore, if this character exists, Siegel's Theorem implies that $d_Q \gg_A (\log Q)^A$ for any $A > 0$. Let $\pi(y; q, a)$ denote the number of primes $p \leqslant y$, $p \equiv a \pmod{q}$. For $n \leqslant 8Q$ not divisible by $d_Q$, $(a, n) = 1$ and $y \geqslant Q_0 := Q^{\log \log Q}$, Theorem 7 of [12] implies

$$\pi(y; n, a) = \frac{y}{\phi(n) \log y}\left(1 + O((\log Q)^{-b})\right)$$

for some absolute constant $b$, $0 < b < 1$. Therefore, for $z \geqslant Q_0$,
$$\sum_{\substack{Q_0 < \varpi \leqslant z \\ \varpi \equiv a \pmod{n}}} \frac{1}{\varpi} = \frac{\log \log z - \log \log Q_0}{\phi(n)} + O\left(\frac{\log \log z}{\phi(n)(\log Q)^b}\right).$$

Write $D = D_1 D_2$, where $D_2$ is the largest odd positive divisor of $D$ (i.e. $D_1 \in \{-2, -1, 1, 2\}$). Since $D_2 \geqslant 3$, quadratic reciprocity implies
$$\left(\frac{D}{\varpi}\right) = (-1)^{(\varpi-1)(D_2-1)/4} \left(\frac{D_1}{\varpi}\right)\left(\frac{\varpi}{D_2}\right) = i_\varpi \left(\frac{\varpi}{D_2}\right),$$
where $i_\varpi \in \{-1, 1\}$ depends only on $\varpi$ modulo 8. Therefore, for numbers $i_e \in \{-1, 1\}$,
$$\sum_{Q_0 < \varpi \leqslant z} \frac{1}{\varpi}\left(\frac{D}{\varpi}\right) = \sum_{e=1,3,5,7} i_e \sum_{\substack{1 \leqslant a \leqslant 8D_2 \\ a \equiv e \pmod 8}} \left(\frac{a}{D_2}\right) \sum_{\substack{Q_0 < \varpi \leqslant z \\ \varpi \equiv a \pmod{8D_2}}} \frac{1}{\varpi}$$
$$= \frac{\log \log z - \log \log Q_0}{\phi(8D_2)} \sum_e i_e \sum_{a'=1}^{D_2} \left(\frac{8a' + e}{D_2}\right) + O\left(\frac{\log \log z}{(\log Q)^b}\right).$$

Here we have used the fact that $\left(\frac{a}{D_2}\right) = 0$ if $(a, D_2) > 1$. The error term is $O(1)$ by hypothesis, and the sum on $a'$ is zero, since as $a'$ runs from 1 to $D_2$, $8a'$ runs over a complete system of residues modulo $D_2$. The lemma now follows from the fact that $1 + h = \exp\{h + O(h^2)\}$ whenever $-\frac{1}{2} \leqslant h \leqslant \frac{1}{2}$. □

## 3. Proof of Theorem 2, Part I

LEMMA 3.1. *If $A(m) = k$ and $S(m)$ is true, then $A(mm') = k + 2$ for some $m'$.*

*Proof.* Suppose $\phi(\alpha_1) = \cdots = \phi(\alpha_k) = m$. Suppose $x$ is sufficiently large, so that $x > 2\alpha_i$ for every $i$ and there are $\gg x/\log^2 x$ primes $q \in (x/2, x]$ such that $qm + 1$ is prime. For each such $q$, the equation $\phi(n) = q(q-1)m$ has at least $k + 2$ trivial solutions, namely $n = q^2 \alpha_i$ ($i = 1, \ldots, k$), $n = q(qm + 1)$ and trivial $n = 2q(qm + 1)$. Let $B(x)$ denote the number of such $q$ giving rise to additional solutions. If $\phi(n) = q(q-1)m$, then clearly $q^3 \nmid n$. If $q^2 \mid n$, then $\phi(n/q^2) = m$, which implies $n$ is one of the trivial solutions. If $q \mid n$, $q^2 \nmid n$, then $\phi(n/q) = qm$. There is a prime $r \mid (n/q)$ with $q \mid (r-1)$, i.e. $r = dq + 1$ for some $d \mid m$. If $d = m$ then $n$ is one of the trivial solutions $mpq + 1$ or $2(mpq + 1)$. By Lemma 2.2, the number of $q$ with $dq + 1$ for some other $d$ is $O(x(\log x)^{-3})$. The last case, counting solutions with $(n, q) = 1$, is codified in a lemma, since we will need the argument in Part II. It follows from the next lemma that the number of $q$ with $\phi(n) = q(q-1)m$ for some $n$ with $(n, q) = 1$ is $O(x(\log x)^{-2.1})$. Therefore, $B(x) = O(x(\log x)^{-2.1})$ and hence for some $q$, $A(q(q-1)m) = k + 2$. □

LEMMA 3.2. *Suppose $b = 1$ or $b$ is a prime less than $y^{10}$. Let $N(y)$ be the number of primes $q \in (y/2, y]$ satisfying* (i) *$qbm + 1$ is prime, and* (ii) *for some $n$ with $(n, q) = 1$, $\phi(n) = mq(q-1)$. Then*
$$N(y) \ll y(\log y)^{-2.1}.$$

*Proof.* First, we remove from consideration those $q$ for which $q - 1$ has a divisor $s$ with $s \geqslant \exp\{(\log \log y)^2\}$ and $P^+(s) \leqslant s^{1/\log \log y}$. By Lemma 2.4, the number of such $q$ is $O(y/\log^{10} y)$. Consider the equation $\phi(n) = mq(q-1)$. There is a number $r \mid n$ with $q \mid r - 1$. Say $r = 1 + qs_0$ and let $w_0 = \phi(n/r)$. For some factorization of $m$ as $m = m_1 m_2$, we have $m_1 \mid s_0$, $m_2 \mid w_0$. Fix $m_1$ and $m_2$, and let $s = s_0/m_1$, $w = w_0/m_2$, so that $sw = q - 1$. The number of $q$ in question is at most the number of pairs $(s, w)$ for which $y/2 - 1 \leqslant sw \leqslant y$, $wm_2$ is a totient (value taken by $\phi(n)$), and the following three numbers are prime:

(3.1) $\quad ws + 1 = q, \qquad (ws + 1)mb + 1 = qmb + 1, \qquad (ws + 1)sm_1 + 1 = r.$

For brevity write $E(c) = \exp\{(\log y)^c\}$. For fixed $s$, Lemma 2.2 implies that the number of $w \leqslant y/s$ with each number in (3.1) prime is
$$\ll \frac{y}{s} \frac{(\log \log y)^3}{\log^3(y/s)}.$$

Therefore, the number of pairs $(s, w)$ with $s \leqslant E(0.89)$ is $O(y(\log y)^{-2.1})$. Otherwise $w \leqslant y/E(0.89)$. By Lemma 2.1 and (2.1), for fixed $w \leqslant y/E(0.89)$, the number of $s$ with each number in (3.1) prime is

(3.2) $\qquad \ll \dfrac{y(\log \log y)^3}{w(\log y)^{2.67}} \prod_{\varpi \leqslant z} \left(1 - \dfrac{1}{\varpi}\left(\dfrac{m_1^2 - 4m_1 w}{\varpi}\right)\right),$

where $z = E(0.89)^{1/2}$.

First, suppose that $w \in (U/2, U]$ where $U$ is a power of 2 less than $2E(1/2)$. Let $Q = 8mE(1/2)$ and recall the definition of $d_Q$ (Lemma 2.8). The number of $w$ with $d_Q \mid (m_1^2 - 4m_1 w)$ is clearly
$$\leqslant \frac{4m_1 U}{d_Q} + 1 \ll \frac{U}{(\log Q)^2} \ll \frac{U}{\log z}.$$

For the remaining $w$, Lemma 2.8 implies that
$$\prod_{Q^{\log \log Q} \leqslant \varpi \leqslant z} \left(1 - \frac{1}{\varpi}\left(\frac{m_1^2 - 4m_1 w}{\varpi}\right)\right) \ll 1.$$

Hence, since $w$ is a totient, Lemma 2.3 gives

$$\sum_{U/2 < w \leqslant U} \frac{1}{w} \prod_{\varpi \leqslant z} \left(1 - \frac{1}{\varpi}\left(\frac{m_1^2 - 4m_1 w}{\varpi}\right)\right) \ll \frac{(\log y)^{1/2} \log \log y}{(\log U)^{0.999}}.$$

By (3.2), summing on $U$ gives a total of $O(y(\log y)^{-2.16})$ possible pairs $(s,w)$.

Lastly, consider $E(1/2) \leqslant w \leqslant y/E(0.89)$. The product in (3.2) for $\varpi > E(0.48)$ is $O((\log y)^{0.41})$. By our initial assumption about $q$, $w$ has a prime factor $> E(0.49)$. Thus $n/r$ has a prime divisor $r_2 > E(0.49)$. We can assume that $r_2^2 \nmid n$, for the number of $q$ divisible by a factor $r_2(r_2-1)$ for such large $r_2$ if $O(y/E(0.49))$. Write $w^* = \phi(n/(rr_2))$, so that $w = (r_2-1)w^*$. Fix $w^*$ and put $r_2 \in (U/2, U]$, where $U$ is a power of 2 in $[E(0.49), 2y/(sw^*)]$. By Lemma 2.7,

$$\sum_{U/2 < r_2 \leqslant U} \frac{1}{r_2} \prod_{\varpi \leqslant E(0.48)} \left(1 - \frac{1}{\varpi}\left(\frac{m_1^2 - 4m_1 w^*(r_2-1)}{\varpi}\right)\right) \ll \frac{(\log \log y)^2}{\log U}.$$

Summing on $U$ and again using (3.2) and Lemma 2.3, we see that the total number of pairs $(s,w)$ is $O(y(\log y)^{-2.25})$. Thus, for fixed $m_1, m_2$, the number of possible $q$ is $O(y(\log y)^{-2.1})$. Since the number of choices for $(m_1, m_2)$ is $O(1)$, the lemma follows. $\square$

## 4. Part II, the case $(n, pq) > 1$

We require one more sieve lemma, needed only for a special case of the proof for Part II.

LEMMA 4.1. *Suppose $r, s, t, u$ are positive integers, and $x, y, z$ are positive real numbers satisfying $x \geqslant y \geqslant z \geqslant 2$ and $\log(rstu) \leqslant \log^3 x$. Let*

$$w_{g,h} = (rg+1)(sh+1)$$

*and*

$$F(g,h) = (tw_{g,h}+1)(ughw_{g,h}+1)ghw_{g,h}.$$

*Then*

$$\#\{x < g \leqslant 2x, y < h \leqslant 2y : P^-(F(g,h)) > z\} \ll xy\frac{(\log \log x)^6}{(\log z)^6}.$$

*Proof.* Let $\varrho(\varpi)$ denote the number of solutions of $F(g,h) \equiv 0 \pmod{\varpi}$. If $\varrho(\varpi) = \varpi^2$ for some $\varpi \leqslant z$ then the left side is zero. Otherwise by the sieve fundamental lemma (Theorem 2.5 of [13]),

$$\#\{x < g \leqslant 2x, y < h \leqslant 2y : P^-(F(g,h)) > z\} \ll xy \prod_{\varpi \leqslant z}\left(1 - \frac{\varrho(\varpi)}{\varpi^2}\right).$$

We next show that for $\varpi \nmid 30rstu$,

(4.1) $$-12 - 2\varpi^{3/4} \leqslant \varrho(\varpi) - 6\varpi \leqslant -10 + 2\varpi^{3/4}.$$

From (4.1), the lemma follows at once, since then

$$\prod_{\varpi \leqslant z}\left(1 - \frac{\varrho(\varpi)}{\varpi^2}\right) \ll \left(\frac{30rstu}{\phi(30rstu)}\right)^6 \prod_{\varpi \leqslant z}(1 - 1/\varpi)^6$$
$$\ll \frac{(\log\log x)^6}{(\log z)^6}.$$

To prove (4.1), suppose $\varpi \leqslant z$, $\varpi \nmid 30rstu$. Let $M_1, M_2, M_3, M_4$ denote the set of solutions mod $\varpi$ of $gh \equiv 0$, $w_{g,h} \equiv 0$, $tw_{g,h} + 1 \equiv 0$ and $ughw_{g,h} + 1 \equiv 0$, respectively. It is easily verified that $M_1 \cap M_4 = M_2 \cap M_3 = M_2 \cap M_4 = \emptyset$, $|M_1 \cap M_2| = |M_1 \cap M_3| = 2$ and $|M_3 \cap M_4| \leqslant 2$. Also, $|M_1| = |M_2| = 2\varpi - 1$ and $|M_3| = \varpi - 1$. By inclusion-exclusion,

$$\varrho(\varpi) = |M_1 \cup M_2 \cup M_3 \cup M_4| = 5\varpi - 7 + |M_4| - |M_3 \cap M_4|.$$

Since the number of solutions of $q(qR + 1) \equiv a \pmod{\varpi}$ is $1 + \left(\frac{1+4aR}{\varpi}\right)$, we have

$$|M_4| = \sum_{a=1}^{\varpi-1}\left(1 + \left(\frac{1+4ar}{\varpi}\right)\right)\left(1 + \left(\frac{1 - 4sa^{-1}u^{-1}}{\varpi}\right)\right)$$
$$= \varpi - 3 + \sum_{a=1}^{\varpi-1}\left(\frac{1+4ar}{\varpi}\right)\left(\frac{a}{\varpi}\right)\left(\frac{a - 4su^{-1}}{\varpi}\right) =: \varpi - 3 + E.$$

By Theorem 1 of Davenport [5], we have $|E| \leqslant 2\varpi^{3/4}$ and (4.1) follows. □

We next combine Lemmas 2.7 and 2.8 into an estimate which has enough generality for our applications.

LEMMA 4.2. *Suppose $\delta > 0$, $z \geqslant 100$. Suppose $A, B, Q, R, U$ are positive integers satisfying $(\log z)^\delta \leqslant \log Q \leqslant (\log z)^3$, $\log(ABRU) < \log^3 z$, $U \geqslant 10$, $(AB, d_Q) \leqslant d_Q^{1/2}$ and $Q \geqslant AU(AU + B)$. Suppose either $D(s) = A(A - sB)$ or $D(s) = sA(sA - B)$. Let $\mathscr{S}_1$ denote the set of primes $s \leqslant U$, and $\mathscr{S}_2$ the set of primes $s \leqslant U$ with $sR + 1$ prime. If $1 \leqslant c \leqslant 2$ and $k \in \{1, 2\}$, then*

$$\sum_{s \in \mathscr{S}_k} \prod_{\varpi \leqslant z}\left(1 - \frac{1}{\varpi}\left(\frac{D(s)}{\varpi}\right)\right)^c \ll_\delta \frac{U}{(\log U)^k}\left(\frac{\log Q}{\log U}\right)^c (\log\log z)^{4+4c} + \log^5 z.$$

*Proof.* Partition $\mathscr{S}_i$ into $\mathscr{S}'$, the set of $s \in \mathscr{S}_i$ for which $d_Q \nmid AD(s)$, and $\mathscr{S}''$, the remaining $s$. If $D(s) = A(A - sB)$, then $\mathscr{S}''$ is contained in a single residue class modulo $g$, where $g = d_Q/(d_Q, AB) \geqslant d_Q^{1/2}$. Recall that $16 \nmid d_Q$ and $d_Q$ is not divisible by the square of any odd prime. If $D(s) = sA(sA - B)$, then $(s, d_Q) > 1$ for at most $\omega(d_Q) \leqslant \log Q \leqslant \log^3 z$ primes $s$ (Here we use $d_Q \leqslant 8Q$). For other $s$, $d_Q | A(sA - B)$ for $s$ lying in a single residue class modulo $g = d'/(d', A)$, where $d' = d_Q/(d_Q, A)$. Here $(d', A) \leqslant 8$, therefore $g \geqslant \frac{1}{8} d_Q^{1/2}$. Thus, in every case

$$|\mathscr{S}''| \ll \frac{U}{g} + 1 + \log^3 z \ll \frac{U}{(\log Q)^{8/\delta}} + \log^3 z,$$

so the contribution to the sum in the lemma from $s \in \mathscr{S}''$ is

$$O(U(\log z)^{-6} + \log^5 z) = O(U/\log^2 U + \log^5 z).$$

By Lemma 2.8, if $s \in \mathscr{S}'$ and $z \geqslant Q_0 := Q^{\log \log Q}$ then

$$\prod_{Q_0 < \varpi \leqslant z} \left(1 - \frac{1}{\varpi}\left(\frac{D(s)}{\varpi}\right)\right)^c \ll_\delta 1.$$

Set $y = U^{1/(20(\log \log U)^2)}$. Trivially we have

$$\prod_{y < \varpi \leqslant Q_0} \left(1 - \frac{1}{\varpi}\left(\frac{D(s)}{\varpi}\right)\right)^c \ll \left(\frac{\log Q_0}{\log y}\right)^c \ll \left(\frac{\log Q}{\log U}\right)^c (\log \log z)^{3c}.$$

Therefore, for any admissible $z$, Lemma 2.7 implies

$$\sum_{s \in \mathscr{S}'} \prod_{\varpi \leqslant z} \left(1 - \frac{1}{\varpi}\left(\frac{D(s)}{\varpi}\right)\right)^c$$

$$\ll \left(\frac{\log Q}{\log U}\right)^c (\log \log z)^{3c} \sum_{s \in \mathscr{S}'} \prod_{\varpi \leqslant \min(y,z)} \left(1 - \frac{1}{\varpi}\left(\frac{D(s)}{\varpi}\right)\right)^c$$

$$\ll \frac{U}{(\log U)^k} \left(\frac{\log Q}{\log U}\right)^c (\log \log z)^{4+4c}. \qquad \square$$

As mentioned in the introduction, for most primes $p \leqslant x$, the prime factors of $p - 1$ are nicely distributed. We prove a simple result of this kind below. For brevity write $E(z) = \exp\{(\log x)^z\}$.

LEMMA 4.3. *Let $N \geqslant 10$ be an integer, put $\delta = 1/N$ and suppose $x^{1/10} \leqslant y \leqslant x$. The number of primes $p \leqslant y$ not satisfying*

$$\Omega(p - 1, E(j\delta), E((j+1)\delta)) \leqslant (\delta + \delta^2) \log \log x \qquad (0 \leqslant j \leqslant N - 1)$$

*is*

$$\ll_\delta y(\log x)^{-1-\delta^3/3}.$$

*Proof.* Assume that $x$ is sufficiently large with respect to $\delta$. First, the number of $p \leqslant y$ for which $P^+(p-1) \leqslant y^{1/\log \log x}$ is $O(y/\log^{10} x)$ by Lemma 2.4. For the remaining $p$, let $q = P^+(p-1)$ and $p-1 = qa$. For fixed $a$, the number of $q \leq y/a$ with $qa + 1$ prime is, by Lemma 2.2,

$$\ll \frac{y \log \log x}{a \log^2(y/a)} \ll \frac{y(\log \log x)^3}{a \log^2 x}.$$

For each $j$, let $I_j = [E(j\delta), E((j+1)\delta)]$ and let $B_j$ be the sum of $1/a$ over $a \leq y$ with at least $(\delta + \delta^2) \log \log x - 1$ prime factors in $I_j$. We have

$$B_j \leq (1+\delta)^{1-(\delta+\delta^2)\log\log x} \sum_{a \leq x} \frac{(1+\delta)^{\Omega(a, E(j\delta), E(j\delta+\delta))}}{a}$$

$$\leq 2(\log x)^{-(\delta+\delta^2)\log(1+\delta)} \prod_{\substack{p \leq x \\ p \notin I_j}} \left(\frac{p}{p-1}\right) \prod_{p \in I_j} \left(1 - \frac{1+\delta}{p}\right)^{-1}$$

$$\ll (\log x)^{1+\delta^2-(\delta+\delta^2)\log(1+\delta)} \ll (\log x)^{1-0.4\delta^3},$$

and the lemma follows. $\square$

Suppose that Hypothesis $S(m)$ (following Lemma 1.2) is false. By Lemma 1.2, for an unbounded set of $x$, there are $\gg x/\log^2 x$ primes $s \in (x/2, x]$ for which $s = 1 + pqm$, where $p$ and $q$ are primes and $x^{1/10} \leq p \leq x^{1/2}$. Denote by $\mathscr{A}$ the set of such pairs of primes $(p, q)$. We will show that for $x$ sufficiently large, there is some such pair $(p, q)$ giving $A(mpq(p-1)(q-1)) = k+2$, with

(4.2) $$\phi(n) = mpq(p-1)(q-1)$$

having only the trivial solutions $n = p^2 q^2 \alpha_i$ $(i = 1, \ldots, k)$, $n = pq(mpq + 1)$ and $n = 2pq(mpq + 1)$. Before proceeding, we must cull from $\mathscr{A}$ certain pairs $(p, q)$ possessing undesirable properties. Let $F(p)$ be the number of $q$ with $(p, q) \in \mathscr{A}$ and let $G(q)$ be the number of $p$ such that $(p, q) \in \mathscr{A}$. By Lemma 2.2, we have

(4.3) $$F(p) \ll \frac{x}{p \log^2 x}, \qquad G(q) \ll \frac{x}{q \log^2 x}.$$

By (4.3) and the Brun-Titchmarsh inequality, if $K$ is a large enough constant (depending only on $m$),

$$\#\{(p, q) \in \mathscr{A} : p \equiv 1 \pmod{2^K} \text{ or } q \equiv 1 \pmod{2^K}\} \leq \tfrac{1}{2}|\mathscr{A}|.$$

Denote by $\mathscr{A}'$ the set of remaining pairs. Define $\nu$ by $2^\nu \mid m$, $2^{\nu+1} \nmid m$. If $n$ has $h$ distinct prime factors, then $2^{h-1} \mid \phi(n)$. Consequently, if (4.2) holds and $(p, q) \in \mathscr{A}'$, $\omega(n) \leq C$ where $C = 2K + 1 + \nu$. Let

(4.4) $$N = 100([4\log(2C+2)] + 14 + 2C), \qquad \delta = 1/N.$$

For brevity, write $E_j = \exp\{(\log x)^{1-j\delta}\}$ for $0 \leq j \leq N-1$ and let $E_N = 2$. Define

(4.5) $\qquad Q_j = E_{j-1}^{3\delta \log \log x}, \quad d_j = d_{Q_j}, \qquad (j = 2, \ldots, N).$

Let $\mathscr{B}$ be the set of $(p,q) \in \mathscr{A}'$ satisfying the additional conditions

(4.6) $\qquad d^2 \mid (p-1)(q-1)$ implies $d \leq \log^4 x$,

(4.7) $\qquad dp+1, dq+1, dpq+1$ are composite for each $d \mid m$,

$\qquad\qquad$ except for $mpq+1$,

(4.8) $\qquad p(p-1)+1, q(q-1)+1, pq(p-1)+1, pq(q-1)+1,$

$\qquad\qquad$ are all composite,

(4.9) $\qquad (m(p-1)(q-1), d_j) \leq d_j^{1/2} \quad (j = 2, \ldots, N),$

(4.10) $\qquad p \nmid (q-1),$

and

(4.11) $\qquad \begin{cases} \Omega(p-1, E_{j+1}, E_j) \leq (\delta + \delta^2) \log \log x & (0 \leq j \leq N-1) \\ \Omega(q-1, E_{j+1}, E_j) \leq (\delta + \delta^2) \log \log x & (0 \leq j \leq N-1). \end{cases}$

First, the number of $(p,q)$ failing (4.6) does not exceed

$$\sum_{d > \log^4 x} \sum_{a \mid d^2} \sum_{n < x, a \mid n} \frac{ax}{nd^2} \ll x \log x \sum_{d > \log^4 x} \frac{\tau(d^2)}{d^2} \ll x(\log x)^{-2.9}.$$

Since $m$ is fixed, by Lemma 2.2 and (4.3), the number of pairs failing (4.7) is $O(x/\log^3 x)$. For (4.8), given $q$ the number of $p$ with $pq(q-1)+1$ prime is $O(x(\log \log x)^3/(q \log^3 x))$ by Lemma 2.2. Thus the number of pairs $(p,q)$ with $pq(q-1)+1$ prime is $O(x(\log x)^{-2.99})$. Likewise the number of pairs satisfying $pq(p-1)+1$ prime is $O(x(\log x)^{-2.99})$. Also, by Lemma 2.2 and the prime number theorem for arithmetic progressions, for each $q$ the number of $p$ with $p(p-1)+1$ prime is

$$\ll \frac{x}{q \log^3 x} \prod_{p \leq x^{1/20}} \left(1 - \frac{1}{p}\left(\frac{-3}{p}\right)\right) \ll \frac{x}{q \log^3 x}.$$

Thus the number of pairs $(p,q)$ with $p(p-1)+1$ prime is $O(x/\log^3 x)$, and the same bound holds for the number of pairs with $q(q-1)+1$ prime. To deal with (4.9), for each $j$ let $V_j$ denote the set of numbers $a \leq d_j$ with $(a, d_j) > d_j^{1/4}/m$. Then

$$|V_j| = \sum_{d \mid d_j, d > d_j^{1/4}/m} \phi(n/d) \leq \sum_{d \mid d_j, d > d_j^{1/4}/m} n/d$$

$$= \sum_{d \mid d_j, d < md_j^{3/4}} d < \tau(d_j) m d_j^{3/4} \ll d_j^{4/5}.$$

By Lemma 2.8, $d_j \gg (\log Q_j)^{40/\delta} \geq (\log x)^{40}$ and for $x$ large, $m \leq d_j^{1/100}$ for every $j$. Therefore, the number of pairs $(p,q)$ failing (4.9) is at most

$$\sum_{j=1}^{N-1} \sum_{a \in V_j} \left( \sum_{\substack{x^{1/10} < p \leq x^{1/2} \\ p \equiv 1+a \pmod{d_j}}} \frac{x}{p} + \sum_{\substack{x^{1/2} < q \leq x^{9/10} \\ q \equiv 1+a \pmod{d_j}}} \frac{x}{q} \right) \ll \sum_{j=1}^{N-1} |V_j| \left( \frac{x \log x}{d_j} + 1 \right)$$

$$\ll \frac{x}{(\log x)^3}.$$

The number of pairs failing (4.10) is trivially $O(x^{9/10})$. By (4.3), Lemma 4.3 and partial summation, the number of pairs failing (4.11) is at most

$$\sum_{\text{bad } p} F(p) + \sum_{\text{bad } q} G(q) \ll \frac{x}{(\log x)^{2+\delta^3/3}}.$$

Thus, for $x$ large,

(4.12) $$|\mathscr{B}| \geq \frac{1}{2} |\mathscr{A}'| \gg x/\log^2 x.$$

We will show that for most pairs $(p, q) \in \mathscr{B}$, the equation (4.2) has only the trivial solutions. We break this argument into two parts, the second of which is the most delicate and is proved in the next section.

LEMMA 4.4. *The number of pairs $(p, q) \in \mathscr{B}$ for which (4.2) has a nontrivial solution $n$ with $(n, pq) > 1$ is $O(x(\log x)^{-2.1})$.*

*Proof.* Consider a generic solution $(n, p, q)$ of (4.2) with $(n, pq) > 1$. By the definition of $\mathscr{B}$, $p^3 \nmid n$ and $q^3 \nmid n$. If $(n, p^2q^2) = p^2q^2$, then $\phi(n/p^2q^2) = m$, which has the $k$ trivial solutions $n = p^2q^2\alpha_i$ ($i = 1, \ldots, k$). Next, if $(n, p^2q^2) = p^2q$, then $\phi(n/p^2q) = mq$. This implies that there is some prime $r \mid n$ with $q|(r-1)|mq$, which is impossible by (4.7). Similarly, $(n, p^2q^2) = pq^2$ is impossible. If $(n, p^2q^2) = pq$, then $\phi(n/pq) = mpq$, implying $n$ has a prime divisor of the form $dp + 1$, $dq + 1$ or $dpq + 1$ for some $d \mid m$. By (4.7), $n = pq(mpq + 1)$ or $n = 2pq(mpq + 1)$. If $(n, p^2q^2) = p^2$, then $\phi(n/p^2) = mq(q-1)$. By Lemma 3.2 (with $b = p$), for fixed $p$ the number of possible $q$ is $O((x/p)(\log x)^{-2.1})$; therefore the total number of possible pairs $(p, q)$ is $O(x(\log x)^{-2.1})$. The same argument handles the case $(n, p^2q^2) = q^2$.

Finally, suppose $(n, p^2q^2) = p$ (the symmetric case $(n, p^2q^2) = q$ is dealt with the same way). Then $\phi(n') = mpq(q-1)$, where $n' = n/p$, and there are two cases: (i) there is a prime $r \mid n'$ with $pq \mid (r-1)$; (ii) there are primes $r_1, r_2$ dividing $n'$ with $p \mid (r_1 - 1)$ and $q \mid (r_2 - 1)$.

In case (i) write $r = 1 + pqs_0$ and $w_0 = \phi(n'/r)$, so $s_0 w_0 = m(q-1)$. Fix $m_1$ and $m_2$, where $m = m_1 m_2$, $m_1 \mid s_0$ and $m_2 \mid w_0$. Write $w = w_0/m_2$ and $s = s_0/m_1$, so that $sw = q - 1$. For fixed $s \leq E(3/4)$ and fixed $p$, the number of $w$ with $sw+1 = q$, $psm_1(sw+1)+1 = r$ and $pm(sw+1)+1$ all prime is

$$\ll \frac{x}{sp} \frac{(\log \log x)^3}{\log^3 x}.$$

Summing on $s \leqslant E(3/4)$ and $p$ yields a total of $O(x(\log x)^{-2.24})$ possible pairs $(p,q)$. The other possibility is that $w < q/E(3/4)$. By Lemma 2.1 (with $z = E(0.74)$) and (2.1), for fixed $w$ and $p$, the number of $s$ is

$$\ll \frac{x}{pw} \frac{(\log \log x)^3}{(\log x)^{2.22}} \prod_{\varpi \leqslant E(0.74)} \left(1 - \frac{1}{\varpi}\left(\frac{m_1^2 p^2 - 4m_1 pw}{\varpi}\right)\right).$$

Put $p \in (U/2, U]$, where $U$ is a power of 2 between $x^{1/10}$ and $x^{1/2}$. By Lemma 2.7, the above product has average $O((\log \log x)^2)$ over such $p$. Summing over $U$, it follows that for fixed $w$, the number of pairs $(p, s)$ is $O((x/w)(\log x)^{-2.21})$. Since $wm_2$ is a totient, Lemma 2.3 implies that the number of pairs $(p, q)$ is $O(x(\log x)^{-2.2})$. Therefore there are $O(x(\log x)^{-2.2})$ pairs $(p,q)$ counted in case (i), since there are $O(1)$ choices for $m_1$ and $m_2$.

In case (ii), write

$$r_1 = 1 + ps_1, \quad r_2 = 1 + qs_2, \quad n' \mid r_1 r_2 t, \quad w_0 = \phi(t).$$

Fix $m_1, m_2, m_3$, where $m = m_1 m_2 m_3$, $m_1 \mid s_1$, $m_2 \mid s_2$ and $m_3 \mid w_0$. Let $s_3 = s_1/m_1$, $s_4 = s_2/m_2$ and $w = w_0/m_3$. The number of possible pairs $(p, q)$ is at most the number of quadruples $(s_3, s_4, w, p)$ where $p$ is prime, $wm_3$ is a totient, and the following four numbers are prime:

(4.13)
$$ws_3 s_4 + 1 = q,$$
$$ps_3 m_1 + 1 = r_1,$$
$$(ws_3 s_4 + 1)s_4 m_2 + 1 = r_2,$$
$$(ws_3 s_4 + 1)pm + 1 = pqm + 1.$$

Suppose $U$ is a power of 2, $x^{1/10} \leqslant U \leqslant 2x^{1/2}$ and $U/2 < p \leqslant U$. We will count separately the quadruples $(s_3, s_4, w, p)$ in each of three classes: 1) those with $w \geqslant E(0.99)$; 2) those with $w < E(0.99)$, $s > E(0.8)$; 3) those with $w < E(0.99)$ and $s \leqslant E(0.8)$.

First, suppose $w \geqslant E(0.99)$. By Lemma 2.2, given $(s_3, s_4)$, the number of $p$ with $ps_3 m_1 + 1$ prime is $O(U(\log \log x)^2/\log^2 x)$ and for each $p$ the number of $w$ with the remaining three numbers in (4.13) all prime is

$$\ll \frac{x(\log \log x)^3}{Us_3 s_4 \log^3(x/(Us_3 s_4))} \ll \frac{x(\log \log x)^3}{Us_3 s_4 (\log x)^{2.97}}.$$

By (4.11), $\Omega(s_3 s_4) \leqslant 1.01 \log \log x$, so

$$\sum_{s_3, s_4} \frac{1}{s_3 s_4} \leqslant \sum_{\substack{m \leqslant x \\ \Omega(m) \leqslant 1.01 \log \log x}} \frac{\tau(m)}{m} \ll (\log x)^{1+1.01 \log 2}.$$

Therefore, the number of quadruples in the first class is $O(x(\log x)^{-2.2})$, since there are $O(\log x)$ possibilities for $U$.

Secondly, suppose $w < E(0.99)$ and $s_3 \geqslant E(0.8)$. Applying Lemma 2.2 twice, for each pair $(s_4, w)$ the number of possible $(p, s_3)$ with $p$ and the numbers in (4.13) all prime is

$$\ll \frac{x}{ws_4} \frac{(\log \log x)^4}{(\log^3 x) \log^2(x/(Us_4 w))} \ll \frac{x(\log \log x)^4}{ws_4 (\log x)^{4.6}}.$$

Since $wm_3$ is a totient, by Lemma 2.3 and partial summation we have

$$\sum_w \frac{1}{w} \sum_{s_4} \frac{1}{s_4} \ll (\log x)^{1.01}.$$

Summing over $U$ gives $O(x(\log x)^{-2.58})$ quadruples in the second class.

Finally, suppose $s_3 < E(0.8)$ and $w < E(0.99)$. By Lemma 2.2, for each pair $(s_3, w)$, the number of $(p, s_4)$ with $p$, $ps_3 m_1 + 1 = r_1$, $s_3 s_4 w + 1 = q$ and $(s_3 s_4 w + 1)pm + 1$ all prime is

$$\ll \frac{x(\log \log x)^4}{s_3 w \log^2 x \log^2(x/(Us_3 w))} \ll \frac{x(\log \log x)^4}{s_3 w \log^4 x}.$$

We have ignored the fact that $r_2$ is prime, because $r_2$ is a quadratic function of $s_4$. As in the previous case, $wm_3$ being a totient implies that the sum of $1/w$ is $\ll (\log x)^{0.01}$. Since the sum of $1/s_3$ is $O((\log x)^{0.8})$, there are $O(x(\log x)^{-2.18})$ quadruples in the third class. The number of triples $(m_1, m_2, m_3)$ is $O(1)$, so there are $O(x(\log x)^{-2.18})$ pairs $(p, q)$ counted in case (ii). □

## 5. Part II, the case $(n, pq) = 1$

In this argument we use for the first time the fact that the prime factors of $p-1$ and $q-1$ are nicely distributed (i.e. the full strength of (4.11)). Recall that $E_j = \exp\{(\log x)^{1-j\delta}\}$ for $0 \leqslant j \leqslant N-1$ and $E_N = 2$.

We may assume that any prime divisor of $n$ which is $\geqslant E_{N-1}$ divides $n$ to the first power, for the number of $(p, q)$ with $r(r-1) \mid m(p-1)(q-1)$ for some $r \geqslant E_{N-1}$ is $O(x/E_{N-1})$. Denote by $\mathscr{B}'$ the set of remaining pairs $(p, q)$.

Let $r_1, \ldots, r_M$ be the prime factors of $n$ satisfying

$$P^+(r_i^*) > E_{N-1}, \qquad r_i^* := \frac{r_i - 1}{(r_i - 1, pq)}.$$

Clearly $1 \leqslant M \leqslant C$. Let $w = \phi(n/(r_1 \cdots r_M))$, so that

$$P^+(w^*) \leqslant E_{N-1}, \qquad w^* = \frac{w}{(w, pq)}.$$

In our new notation, we have

(5.1) $$r_1^* \cdots r_M^* w^* = m(p-1)(q-1).$$

We will partition the solutions of (5.1) according to the value of $M$ (the number of "large" prime factors of $n$), the size of the largest prime factor of each $r_i^*$, and the location of the factors $p$ and $q$ among the factors $r_i - 1$ and $w$. For convenience, set $r_0 = w+1$. Define numbers $i'$ and $i''$ by $p \mid (r_{i'}-1)$, $q \mid (r_{i''}-1)$. It is possible that $i' = i''$. Let $I_j = [E_j, E_{j-1})$ and define numbers $j_i$ by $P^+(r_i^*) \in I_{j_i}$. The primes $r_i$ may be ordered so that $1 = j_1 \leqslant j_2 \leqslant \cdots \leqslant j_M$. The fact that $j_1 = 1$ follows from the $j = 0$ case of (4.11).

By (4.8), it follows that $p - 1 \neq r_{i'}^*$ and $q - 1 \neq r_{i''}^*$. Thus we may define

(5.2) $$S' = P^+\left(\frac{(p-1)r_{i'}^*}{(p-1, r_{i'}^*)^2}\right), \qquad S'' = P^+\left(\frac{(q-1)r_{i''}^*}{(q-1, r_{i''}^*)^2}\right).$$

In other words, $S'$ is the largest prime dividing $p-1$ and $r_i^*$ to different powers. Define $j', j''$ by $S' \in I_{j'}$, $S'' \in I_{j''}$, so that $j'$ and $j''$ measure the rough sizes of $S'$ and $S''$. By (4.11), if $j' > 1$ then the prime factors of $p-1$ which are $\geqslant E_1$ also divide $r_{i'}^*$, so that $j_{i'} = 1$. Similarly, if $j'' > 1$ then $j_{i''} = 1$. Furthermore, if either $j' > 1$ or $j'' > 1$ then there are at least two values of $i$ with $j_i = 1$. Thus, without loss of generality, if $j' > 1$ we may assume $i' = 1$ and if $j'' > 1$, we may assume $i'' = 2$.

We count separately the number of $(p,q) \in \mathcal{B}'$ giving solutions of (5.1) with $M$, $j_1, \ldots, j_M$, $i'$, $i''$, $j'$ and $j''$ fixed. The number of choices for these parameters is $O(1)$ since it only depends on $m$.

To reap the full benefit of (4.11), we break each factor in (5.1) into $N$ pieces as follows:

$$m(p-1) = \prod_{j=1}^{N} \pi_j, \quad q-1 = \prod_{j=1}^{N} \theta_j, \quad r_i^* = \prod_{j=1}^{N} \varrho_{i,j} \quad (1 \leqslant i \leqslant M),$$

where $\pi_j$, $\theta_j$ and $\varrho_{i,j}$ are composed only of primes in $I_j$. We can suppose $x$ is large enough so that $P^+(m) < E_{N-1}$, i.e. $m \mid \pi_N$. Some of these variables may equal 1. In particular, for each $i$ we have $\varrho_{i,j} = 1$ for $j < j_i$. By the definition of $j'$ and $j''$ we also have

(5.3) $$\pi_j = \varrho_{1,j} \quad (j < j'), \qquad \theta_j = \varrho_{2,j} \quad (j < j'').$$

Let

(5.4) $$a_j = |\{i \geqslant 1 : j_i \leqslant j\}| \qquad (j = 1, \ldots, N-1),$$
$$b_0 = |\{i \geqslant 1 : (r_i - 1, pq) > 1\}|,$$
$$b_j = |\{i \geqslant 1 : j_i = j, (r_i - 1, pq) = 1\}|, \qquad (j = 1, \ldots, N-1).$$

For each $j$, $a_j$ is the maximum number of $\varrho_{i,j}$ which can be $> 1$. The estimation of the number of solutions of (5.1) proceeds in $N$ steps; at step $j$ we count the possibilities for $\pi_j, \theta_j$ and the $\varrho_{i,j}$, assuming the variables consisting of smaller primes have been fixed. The possibility that $j' > 1$ or $j'' > 1$ will force us to express some bounds in terms of products involving the Legendre symbol, and these bounds must be carried over many steps. For convenience, for each $j \geqslant 1$ set

$$(5.5) \qquad R_{i,j} = \prod_{h \geqslant j} \varrho_{i,h}, \quad P_j = \prod_{h \geqslant j} \pi_h, \quad T_j = \prod_{h \geqslant j} \theta_h.$$

If $j' > 1$ define

$$(5.6) \qquad D_1 = R_{1,j'}(R_{1,j'} - 4P_{j'}/m), \quad H_1 = \prod_{\varpi \leqslant E_1} \left(1 - \frac{1}{\varpi}\left(\frac{D_1}{\varpi}\right)\right),$$

and if $j'' > 1$ define

$$(5.7) \qquad D_2 = R_{2,j''}(R_{2,j''} - 4T_{j''}), \quad H_2 = \prod_{\varpi \leqslant E_1} \left(1 - \frac{1}{\varpi}\left(\frac{D_2}{\varpi}\right)\right).$$

Note that by (5.3),

$$(5.8) \qquad \left(\frac{D_1}{\varpi}\right) = \left(\frac{r_1^*(r_1^* - 4(p-1))}{\varpi}\right), \quad \left(\frac{D_2}{\varpi}\right) = \left(\frac{r_2^*(r_2^* - 4(q-1))}{\varpi}\right).$$

Below are the estimates that we will prove:

*Step* 1. Let $w^*$, $P_2$, $T_2$ and $R_{1,2}, \ldots, R_{M,2}$ be given. The number of possibilities for $\pi_1, \theta_1$ and $\varrho_{1,1}, \ldots, \varrho_{M,1}$ is

$$\ll \frac{x}{P_2 T_2} (\log \log x)^{3+b_0+b_1+2a_1} (\log x)^{-4-b_0-b_1+\delta(3+b_0+b_1+2a_1)} X_1 X_2,$$

where

$$X_1 = \begin{cases} 1 & j' = 1 \\ H_1 & j' > 1 \end{cases}, \qquad X_2 = \begin{cases} 1 & j'' = 1 \\ H_2 & j'' > 1 \end{cases}.$$

Next, define

$$Y_j = \begin{cases} 1 & j \neq j' \\ H_1 & j = j' \end{cases}, \qquad Z_j = \begin{cases} 1 & j \neq j'' \\ H_2 & j = j'' \end{cases}.$$

*Step* $j$, $2 \leqslant j \leqslant N - 1$. Let $w^*$, $P_{j+1}$, $T_{j+1}$ and $R_{1,j+1}, \ldots, R_{M,j+1}$ be given. Then

$$\sum \frac{Y_j Z_j}{\pi_j \theta_j} \ll (\log \log x)^{2b_j + 14\zeta_j} (\log x)^{2\delta(1+\delta)\log a_j + 2\delta + b_j(j\delta - 1) + 4\zeta_j \delta},$$

where the sum is over all possible $\pi_j, \theta_j$ and $\varrho_{1,j}, \ldots, \varrho_{M,j}$, and $\zeta_j = 1$ if $j = j'$ or $j = j''$ and $\zeta_j = 0$ otherwise.

*Step N.* We have
$$\sum \frac{Y_N Z_N}{\pi_N \theta_N} \ll (\log \log x)^{14} (\log x)^{5\delta + 2\delta(1+\delta) \log(2M+2)},$$
where the sum is over $w^*$, $\pi_N$, $\theta_N$, and $\varrho_{i,N}$ $(1 \leqslant i \leqslant M)$.

Using the fact that $b_0 + \cdots + b_{N-1} = M$ and combining all $N$ steps, the total number of solutions of (5.1) is
$$\ll x (\log \log x)^{2M+45+a_1} (\log x)^Y,$$
where
$$Y = -4 - M + 2\delta^2 \left( \sum_{j=2}^{N-1} \log a_j + \log(2M+2) \right) + \delta X,$$
$$X = 16 + b_0 + b_1 + 2a_1 + \sum_{j=2}^{N-1} (jb_j + 2 + 2\log a_j) + 2\log(2M+2).$$

Since $a_j \leqslant M$ for each $j$, the coefficient of $\delta^2$ is at most $2N \log(2M+2)$. Also, $b_0 \leqslant 2$; thus
$$X \leqslant 2\log(2M+2) + 14 + 2N + 2M + \sum_{j=1}^{N-1} (2\log a_j + jb_j).$$

The number of $a_j = i$ is $j_{i+1} - j_i$. Also, $j_1 = 1$ and $2\log 2 > 1$; thus, if $M \geqslant 2$ then
$$\sum_{j=1}^{N-1} (2\log a_j + jb_j) \leqslant \sum_{i=1}^{M} j_i + 2 \sum_{i=1}^{M-1} (j_{i+1} - j_i) \log i + 2(N - 1 - j_M) \log M$$
$$= 2(N-1) \log M + j_1 + \sum_{i=2}^{M} \left( 1 - 2\log\left(\tfrac{i}{i-1}\right) \right) j_i$$
$$\leqslant 2N \log M + N \sum_{i=3}^{M} \left( 1 - 2\log\left(\tfrac{i}{i-1}\right) \right)$$
$$= N(M - 2 + 2\log 2).$$

The above bound also holds when $M = 1$, since $a_i = 1$ for all $i$ and $b_j = 0$ for $j \geqslant 1$. Therefore, by (4.4), we obtain
$$Y \leqslant -4 + 2\log 2 + \delta(14 + 2M + 4\log(2M+2)) \leqslant -2.603$$
for every possible $M$, $i'$, $i''$, $j'$, $j''$ and $j_1, \ldots, j_M$. In conclusion, the number of solutions of (4.2) with $(n, pq) = 1$ is $O(x(\log x)^{-2.6})$. Together with (4.12) and Lemma 4.4, this completes the proof of Part II of Theorem 2.

We now establish the bounds claimed for each of the $N$ steps.

*Proof of Step* 1 *bound.* Let $I = \{i : j_i = 1\} = \{1, 2, \ldots, a_1\}$. Recall (5.5) and for $i \in I$ write

(5.9) $$\varrho_{i,1} = A_i B_i, \qquad A_i = (\varrho_{i,1}, \pi_1), \qquad B_i = (\varrho_{i,1}, \theta_1).$$

Then

(5.10) $$\pi_1 = \prod_{i \in I} A_i, \qquad \theta_1 = \prod_{i \in I} B_i$$

and

(5.11) $$P^-\left(\prod_{i \in I} A_i B_i\right) \geqslant E_1.$$

We need to count the number of $(A_1, \ldots, B_{a_1})$ satisfying (5.11) such that the following $3 + b_0 + b_1$ numbers are prime:

(5.12) $$\begin{cases} \pi_1 P_1/m + 1 = p, \\ \theta_1 T_1 + 1 = q, \\ (\pi_1 P_1/m + 1)(\theta_1 T_1 + 1)m + 1 = pqm + 1, \\ A_i B_i R_{1,i}(r_i - 1, pq) + 1 = r_i, & (i \in I), \\ R_{1,i}(r_i - 1, pq) + 1 = r_i, & (i = i' \notin I \text{ or } i = i'' \notin I). \end{cases}$$

Recall that $p \mid (r_{i'} - 1)$ and $q \mid (r_{i''} - 1)$. All variables $A_i, B_i$ occur linearly on the left of the above expressions except in a few cases. If $i' \in I$ (i.e. $j_{i'} = 1$), then $A_{i'}$ occurs quadratically in the expression for $r_{i'}$ (fourth line in (5.12)), and if $i'' \in I$ then $B_{i''}$ occurs quadratically in the expression for $r_{i''}$ (fourth line of (5.12)). There are $O(1)$ possibilities for $I_1 = \{i : A_i > 1\}$, $I_2 = \{i : B_i > 1\}$. Fix $I_1, I_2$ and note that both $I_1$ and $I_2$ are nonempty and $I_1 \cup I_2 = I$. Let $J = |I_1| + |I_2|$.

Place each variable which is $> 1$ into a dyadic interval, e.g. $U_i < A_i \leqslant 2U_i$ for $i \in I_1$, $V_i < B_i \leqslant 2V_i$ for $i \in I_2$, where each $U_i, V_i$ is a power of 2. For some ordering of the $J$ variables, we use a sieve to bound the possibilities for the first variable with the other variables undetermined. Then with the first variable fixed, we use the sieve to bound the possibilities for the second variable, with the remaining $J - 2$ variables undetermined, and so on. In most cases the variables may be ordered so as to apply Lemma 2.2 with each. The exceptions are when either $A_{i'}$ appears quadratically in the expression for $r_{i'}$ and is isolated (none of the other $J - 1$ variables appears in the expression for $r_{i'}$), or when $B_{i''}$ has a similar property, or when $i' = i'' = 1$ and $A_{i'}$ and $B_{i''}$ are isolated in the expression for $r_{i'}$. We therefore distinguish three cases: (i) $a_1 = i' = i'' = 1$; (ii) $j' = j'' = 1$ other than case (i); (iii) $j' + j'' > 2$.

In case (i), we have $J = 2$, $j' = j'' = 1$, $b_0 = 1$, $b_1 = 0$, $\pi_1 = A_1$ and $\theta_1 = B_1$. By Lemma 4.1, the number of pairs $A_1, B_1$ satisfying (5.11) such that $p = A_1(P_1/m) + 1$, $q = B_1 T_1 + 1$, $pqm + 1$ and $A_1 B_1 R_{1,1} pq + 1 = r_1$ are prime is $O(U_1 V_1 (\log \log x)^6 (\log E_1)^{-6})$.

In case (ii), we first choose the numbers $A_{i'}$ (if $i' > 0$) and $B_{i''}$ (if $i'' > 0$). Since $j' = j'' = 1$, none of the primes in (5.12) are determined solely by these choices, so no sieving is required. The remaining variables $A_i$ and $B_i$ occur linearly on the left of the expressions in (5.12), so we may apply Lemma 2.2 repeatedly. The number of solutions counted is

$$\ll \left( \prod_{i \in I_1} U_i \prod_{i \in I_2} V_i \right) (\log \log x)^{3+b_0+b_1+J} (\log E_1)^{-(3+b_0+b_1+J)}.$$

For case (iii), suppose $j' > 1$, $i' = 1$, so that $\pi_1 = \varrho_{1,1} = A_1$ and $r_1 = A_1 R_{1,1}(A_1 P_1/m + 1) + 1$. By Lemma 2.1 and (5.8), the number of possible $A_1$ is $\ll U_1 H_1 (\log \log x)^2 (\log E_1)^{-2}$. Similarly, if $j'' > 1$ the number of possible $B_2$ is $\ll U_2 H_2 (\log \log x)^2 (\log E_1)^{-2}$. Lemma 2.2 is applied repeatedly to the other variables yielding a total of

$$\ll \left( \prod_{i \in I_1} U_i \prod_{i \in I_2} V_i \right) (\log \log x)^{3+b_0+b_1+J} (\log E_1)^{-(3+b_0+b_1+J)} X_1 X_2$$

solutions counted.

In every case, with $I_1$ and $I_2$ fixed,

$$\sum_{U_i, V_i} \prod_{i \in I_1} U_i \prod_{i \in I_2} V_i \ll \frac{x}{P_2 T_2} (\log x)^{J-1}$$

and the bound in Step 1 follows.

*Proof of bound for step $j \in [2, N-1]$.* All $j$. There are $b_j$ numbers $i$ with $j_i = j$ and $(r_i - 1, pq) = 1$. Recall (5.2) and let $I$ denote the set of such $i$ with the additional property that $S' \nmid \varrho_{i,j}$ and $S'' \nmid \varrho_{i,j}$. For each $i \in I$, let $W_i = P^+(\varrho_{i,j})$. Fix the partition of these $i$ into $I_1$, the set of $i$ with $W_i \mid \theta_j$, and $I_2$, the set of $i$ with $W_i \mid \pi_j$. There are $O(1)$ such partitions. Also fix

$$\tau := \{\varrho^*_{1,j}, \ldots, \varrho^*_{M,j}; \pi^*_j, \theta^*_j\}.$$

Here, if $j \neq j'$ and $j \neq j''$, we take $\varrho^*_{i,j} = \varrho_{i,j}/W_i$ if $i \in I$ and $\varrho^*_{i,j} = \varrho_{i,j}$ otherwise, $\pi^*_j = \pi_j / \prod_{i \in I_2} W_i$ and $\theta^*_j = \theta_j / \prod_{i \in I_1} W_i$. If $j = j'$, we also divide out the factor $S'$ from either $\pi_j$ or $\theta_j$, as well as the appropriate $\varrho_{i,j}$. If $j = j''$ we divide out the factor $S''$ in the same manner. Let

(5.13) $$t = \varrho^*_{1,j} \cdots \varrho^*_{M,j} = \pi^*_j \theta^*_j$$

and note that by (4.11),

(5.14) $$\Omega(t) \leqslant 2\delta(1+\delta) \log \log x.$$

Our procedure will be to fix $\tau$ and the $W_i$ ($i \in I$), and sum $Y_j Z_j / S' S''$ over $S', S''$ if necessary. Then for each $i \in I$, $W_i$ and $W_i \varrho_{i,j}^* R_{i,j+1} + 1$ are prime, so by Lemma 2.2 and partial summation,

$$(5.15) \qquad \sum_{W_i} \frac{1}{W_i} \ll (\log \log x)^2 (\log x)^{j\delta - 1}.$$

For fixed $t$, let $h(t)$ be the number of $(M+2)$-tuples $\tau$ with product $t$. At most $a_j$ of the numbers $\varrho_{i,j}$ are $> 1$, thus $h(t)$ is at most the number of dual factorizations of $t$ into 2 and $a_j$ factors, respectively. Actually if $j < \max(j', j'')$ we can do better, but this bound suffices. Therefore $h(t) \leqslant (2a_j)^{\Omega(t)}$. Since $2 \notin I_j$, by (5.14) we have

$$(5.16) \qquad \sum_t \frac{h(t)}{t} \leqslant a_j^{2\delta(1+\delta) \log \log x} \sum_t \frac{2^{\Omega(t)}}{t}$$
$$\leqslant (\log x)^{2\delta(1+\delta) \log a_j} \prod_{p \in I_j} (1 - 2/p)^{-1}$$
$$\ll (\log x)^{2\delta(1+\delta) \log a_j + 2\delta}.$$

All that remains is to take care of the sum on $S', S''$.

*Proof of bound for step $j \in [2, N-1]$; $j \neq j'$, $j \neq j''$.* Here there is no sum of $S'$ or $S''$, so we are done by (5.15) and (5.16).

*Proof of bound in step $j \in [2, N-1]$; $j = j' \neq j''$.* Since $j' > 0$, we have $i' = 1$. There are three cases to consider: (i) $S' \mid \varrho_{1,j}$, $S' \nmid \pi_j$; (ii) $S' \mid \pi_j$, $S' \nmid \varrho_{1,j}$, $S' \mid \varrho_{i,j}$ with $j_i < j$; (iii) $S' \mid \pi_j$, $S' \nmid \varrho_{1,j}$, $S' \mid \varrho_{i,j}$ with $j_i = j$; for all cases suppose $S' \in (U, 2U]$ where $U$ is a power of 2 in $I_j$. In case (i), (5.6) gives

$$D_1 = AS'(AS' - 4P_{j'}/m),$$

with $A = R_{1,j'}/S'$ and $P_{j'}/m$ fixed. Therefore, by (4.5), (4.9) and Lemma 4.2 (with $k=1$, $c=1$, $Q = Q_{j'}$),

$$\sum_U \sum_{U < S' \leqslant 2U} \frac{H_1}{S'} \ll \sum_U \frac{\log Q_{j'} (\log \log z)^4}{\log^2 U} \ll (\log \log x)^5 (\log x)^\delta.$$

In case (ii), (5.6) gives

$$D_1 = R_{1,j'}(R_{1,j'} - S'B),$$

with $B = 4P_j/(mS')$ and $R_{1,j'}$ fixed. As in case (i), by Lemma 4.2 (with $k=1$, $c=1$, $Q = Q_{j'}$), we have

$$\sum_U \sum_{U < S' \leqslant 2U} \frac{H_1}{S'} \ll (\log \log x)^5 (\log x)^\delta.$$

For case (iii), we also have
$$D_1 = R_{1,j'}(R_{1,j'} - S'B),$$
with $B = 4P_j/(mS')$ and $R_{1,j'}$ fixed. Here we must sum over $S'$ for which $S'R + 1$ is prime, where $R = R_{i,j}/S'$ is fixed. By (4.5) and Lemma 4.2 (with $k = 2$, $c = 1$), we obtain
$$\sum_U \sum_{U < S' \leqslant 2U} \frac{H_1}{S'} \ll \sum_U \frac{\log Q_{j'} (\log \log x)^4}{\log^3 U} \ll (\log \log x)^5 (\log x)^{\delta - 1 + j\delta}.$$
The desired bounds now follow from (5.15) and (5.16).

A symmetric argument using (5.7) proves the step $j$ bound when $j = j'' \neq j'$.

*Proof of bound in step $j$; $j = j' = j''$; $S' \neq S''$.* By (5.3), $i_{j-1} = 2$. Define $i^*$, $i^{**}$ by $S' \mid \varrho_{i^*,j}$, $S'' \mid \varrho_{i^{**},j}$. Put $S' \in (U', 2U']$, $S'' \in (U'', 2U'']$ where $U'$ and $U''$ are powers of 2 in $[\frac{1}{2}E_j, E_{j-1})$. By the Cauchy-Schwarz inequality,
$$\sum_{U',U''} \sum_{S',S''} \frac{H_1 H_2}{S'S''} \leqslant \sum_{U',U''} \frac{1}{U'U''} \left( \sum_{S',S''} H_1^2 \right)^{1/2} \left( \sum_{S',S''} H_2^2 \right)^{1/2}.$$
For each double sum over $S', S''$, we first fix $S''$ and sum over $S'$, then sum on $S''$, bounding each sum with Lemma 4.2, Lemma 2.2 or the prime number theorem, as appropriate. For example, for the sum on $S'$ in the first sum we use Lemma 4.2 (with $c = 2$) if $D_1$ depends on $S'$. If $i^* \geqslant 3$ then $S'R' + 1$ is also prime, where $R' = R_{i^*,j}/S'$ is fixed, so we use the $k = 2$ case of the lemma. Otherwise we use the $k = 1$ case. If $D_1$ does not depend on $S'$, we use Lemma 2.2 if $i^* \geqslant 3$ and the prime number theorem otherwise. For the sum on $S''$, we use Lemma 4.2 only if $D_1$ depends on $S''$ but not on $S'$. If $i^{**} \geqslant 3$ and $i^{**} \neq i^*$, then we use the $k = 2$ case of Lemma 4.2, since $S''R'' + 1$ is prime, where $R'' = R_{i^{**},j}/S''$ is fixed. Otherwise use the $k = 1$ case. If $D_1$ does not depend on $S''$, or it depends on both $S'$ and $S''$, we apply either Lemma 2.2, if $i^{**} \geqslant 3$ and $i^{**} \neq i^*$, or the prime number theorem, otherwise. In every case we obtain the bound
$$\sum_{S',S''} H_1^2 \ll \frac{U'U'' \log^4 Q_j}{(\log U')^{k_1+3}(\log U'')^{k_2+3}} (\log \log x)^{10},$$
where $k_1 = 1$ if $i^{**} \geqslant 3$ and $i^{**} \neq i^*$, $k_1 = 0$ otherwise, and $k_2 = 1$ if $i^* \geqslant 3$ and $k_2 = 0$ otherwise. The same bound holds for the sum on $H_2^2$, so
$$\sum_{U',U''} \sum_{S',S''} \frac{H_1 H_2}{S'S''} \ll (\log \log x)^{14} (\log x)^{4\delta - (k_1+k_2)(1-j\delta)}.$$
The desired bound again follows at once from (5.15) and (5.16).

*Proof of bound for step* $j \in [2, N-1]$; $j = j' = j''$, $S' = S''$. Since $(\pi_j, \theta_j) = (\varrho_{1,j}, \varrho_{2,j}) = 1$, either $S' \mid \pi_j, S' \mid \varrho_{2,j}$ or $S' \mid \theta_j, S' \mid \varrho_{1,j}$. In either case, by (5.6) and (5.7), $D_1$ and $D_2$ both depend on $S'$ (recall $\tau$ and the $W_i$ for $i \in I$ are fixed). As before we put $S' \in (U, 2U]$ where $U$ is a power of 2 in $[\frac{1}{2}E_{j-1}, E_j)$. By the Cauchy-Schwarz inequality and Lemma 4.2 (with $c = 2$, $k = 1$) we obtain

$$\sum_U \sum_{U < S' \leqslant 2U} \frac{H_1 H_2}{S'} \leqslant \sum_U \frac{1}{U} \left(\sum_{S'} H_1^2\right)^{1/2} \left(\sum_{S'} H_2^2\right)^{1/2}$$
$$\ll (\log \log x)^{14} (\log x)^{4\delta}.$$

Combined with (5.15) and (5.16) gives the desired bound. This concludes the proof of the bounds for steps $j = 2$ through $j = N - 1$.

*Proof of the bound for step* $N$. Let

$$\tau = \{\pi_N^*, \theta_N^*; \varrho_{1,N}^*, \cdots, \varrho_{M,N}^*, w^*\},$$

where the starred variables (except $w^*$) are the quotients of the corresponding unstarred variables and $S'$, $S''$, as appropriate. As before let $t = \pi_N^* \theta_N^*$. We proceed as in steps 2 through $N-1$, except now there are no variables $W_i$ to sum over. The sum over $S'$ and $S''$ is handled exactly as in steps 2 through $N-1$, but the sum over $\tau$ of $1/t$ is handled a bit differently. Since $2 \in I_N$, the procedure in (5.16) will not work. By (4.11), we have $\Omega(t) \leqslant 2\delta(1+\delta) \log \log x + \Omega(m)$. As in the other steps, given $t$, the number of dual factorizations of $t$ into $M+1$ factors and 2 factors is at most $(2M+2)^{\Omega(t)}$. We thus have

$$\sum_\tau \frac{1}{t} \ll (2M+2)^{2\delta(1+\delta) \log \log x} \sum_{P^+(t) \leqslant E_{N-1}} \frac{1}{t}$$
$$\ll (\log x)^{\delta + 2\delta(1+\delta) \log(2M+2)},$$

and the bound for step $N$ follows.

*Acknowledgements.* The author thanks Sergei Konyagin for helpful discussions and the referee for helpful criticisms of the exposition.

University of Texas at Austin, Austin, TX
*Current address*: University of South Carolina, Columbia, SC
*E-mail address*: ford@math.sc.edu

### References

[1] R. D. Carmichael, On Euler's $\phi$-function, Bull. A.M.S. **13** (1907), 241–243.
[2] ———, Note on Euler's $\phi$-function, Bull. A.M.S. **28** (1922), 109–110.
[3] J. R. Chen, On the representation of a large even integer as the sum of a prime and the product of at most two primes, Kexue Tongbao **17** (1966), 385–386.

[4] J. R. Chen, On the representation of a large even integer as the sum of a prime and the product of at most two primes, Sci. Sinica **16** (1973), 157–176.

[5] H. Davenport, On the distribution of quadratic residues (mod $p$), J. London Math. Soc. **6** (1931), 49–54.

[6] ———, *Multiplicative Number Theory*, 2nd ed., Grad. Texts in Math. **74**, Springer-Verlag, New York, 1980.

[7] L. E. Dickson, A new extension of Dirichlet's theorem on prime numbers, Messenger of Math. **33** (1904), 155–161.

[8] P. Erdős, On the normal number of prime factors of $p-1$ and some related problems concerning Euler's $\phi$-function, Quart. J. Math. Oxford **6** (1935), 205–213.

[9] ———, Some remarks on Euler's $\phi$-function, Acta Arith. **4** (1958), 10–19.

[10] K. Ford, The distribution of totients, Ramanujan J. **2** (1998), 67–151.

[11] K. Ford and S. Konyagin, On two conjectures of Sierpiński concerning the arithmetic functions $\sigma$ and $\phi$, in *Number Theory in Progress*, Proc. Internat. Conf. on Number Theory in honor of 60[th] birthday of Andrzej Schinzel (Zakopane, Poland, 1997), Walter de Gruyter, New York, 1999, 795–803.

[12] P. X. Gallagher, A large sieve density estimate near $\sigma = 1$, Invent. Math. **11** (1970), 329–339.

[13] H. Halberstam and H.-E. Richert, *Sieve Methods*, Academic Press, London, 1974.

[14] G.H. Hardy and S. Ramanujan, The normal number of prime factors of a number $n$, Quart. J. Math. **48** (1917), 76–92.

[15] A. Hildebrand and G. Tenenbaum, Integers without large prime factors, J. Théor. Nombres Bordeaux **5** (1993), 411–484.

[16] A. Schinzel, Sur l'equation $\phi(x) = m$, Elem. Math. **11** (1956), 75–78.

[17] ———, Remarks on the paper "Sur certaines hypothèses concernant les nombres premiers", Acta Arith. **7** (1961/62), 1–8.

[18] A. Schinzel and W. Sierpiński, Sur certaines hypothèses concernant les nombres premiers, Acta Arith. **4** (1958), 185–208.